\documentclass[11pt]{amsart}

 \usepackage{graphicx}
 \usepackage[all]{xy}

\begin{document}

\theoremstyle{plain}
    \newtheorem{theorem}{Theorem}

\title[Semilinear clannish algebras from triangulations]{Semilinear clannish algebras associated to triangulations of surfaces with orbifold points: Oberwolfach talk, February 2023}
\author{Daniel Labardini-Fragoso}
\address{Daniel Labardini-Fragoso\newline
Instituto de Matem\'aticas, UNAM, Mexico}
\email{labardini@im.unam.mx}

\begin{abstract}

\noindent This is an extended abstract of my talk at the Oberwolfach Workshop ``Representation Theory of Quivers and Finite-Dimensional Algebras'' (February 12 - February 18, 2023). It is based on a joint work with R. Bennett-Tennenhaus (arXiv:2303.05326).
\end{abstract}

\maketitle


The main aim of this talk is to present \cite{DLF_BTLF}. We construct semilinear clannish algebras for the colored triangulations of a surface with marked points and orbifold points, and prove that they are Morita-equivalent to the Jacobian algebras of the species with potential constructed by Geuenich and myself a few years ago in \cite{DLF_GLF1,DLF_GLF2}.

\vspace{2mm}

\noindent\textbf{Surfaces with marked points and orbifold points.} A \emph{surface with marked points and orbifold points} is a triple $\mathbf{\Sigma}=(\Sigma,\mathbb{M},\mathbb{O})$ consisting of
\begin{itemize}
\item A compact, connected, oriented, two-dimensional real manifold $\Sigma$ with (possibly empty) boundary $\partial\Sigma$;
\item a finite set of \emph{marked points} $\emptyset\neq \mathbb{M}\subseteq\Sigma$ with at least one point from each connected component of $\partial\Sigma$; points in $\mathbb{P}:=\mathbb{M}\setminus\partial\Sigma$ are called \emph{punctures};
\item a (possibly empty) finite set of \emph{orbifold points} $\mathbb{O}\subseteq\Sigma\setminus\partial\Sigma$.
\end{itemize}

An \emph{arc on $\boldsymbol{\Sigma}$} is a curve $k$ on $\Sigma$ that connects either a pair of points in $\mathbb{M}$, or a point in $\mathbb{M}$ and a point in $\mathbb{O}$, and satisfies the following conditions:
\begin{itemize}
\item except for its endpoints, $k$ is disjoint from $\partial\Sigma\cup\mathbb{M}\cup\mathbb{O}$;
\item except possibly for its endpoints, $k$ does not cross itself;
\item $k$ is not homotopically trivial in $\Sigma\setminus(\mathbb{M}\cup\mathbb{O})$ rel $\mathbb{M}\cup\mathbb{O}$;
\item $k$ is not homotopic in $\Sigma$ rel $\mathbb{M}\cup\mathbb{O}$ to a boundary segment of $\Sigma$;
\item $k$ is not a loop closely enclosing a single orbifold point.
\end{itemize} 
There are two types of arcs: those connecting points in $\mathbb{M}$, called \emph{non-pending arcs}, and those connecting a point in $\mathbb{M}$ to a point in $\mathbb{O}$, called \emph{pending arcs}.

Arcs are considered up to homotopy rel $\mathbb{M}\cup\mathbb{O}$.
Two arcs are \emph{compatible} if there are representatives in their homotopy classes rel $\mathbb{M}\cup\mathbb{O}$ that do not intersect in $\Sigma\setminus\partial\Sigma$. A \emph{triangulation} of $\boldsymbol{\Sigma}$ is a maximal collection of pairwise compatible arcs. Each triangulation $\tau$ of $\mathbf{\Sigma}$ splits $\Sigma$ into finitely many triangles. 

\begin{figure}[!ht]
                \caption{Left: The tree types of triangles of a triangulation.
                Right: The cells of the CW-complex $X(\tau)=(X_n(\tau))_{n=0,1,2}$.}\label{DLF_Fig_unpunct_puzzle_pieces}
                \includegraphics[scale=.1]{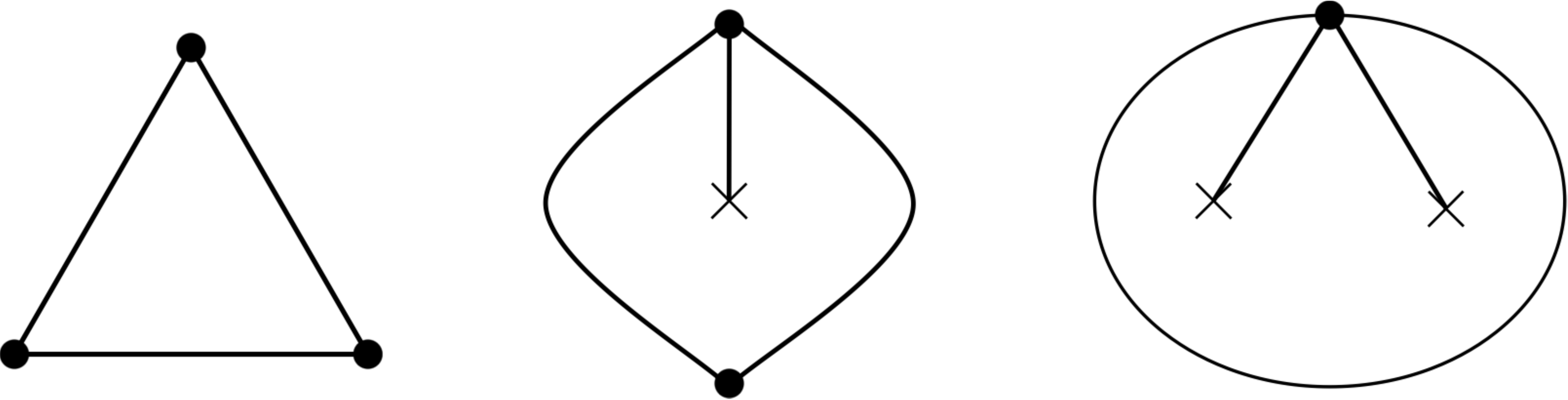}\quad
                 \includegraphics[scale=.11]{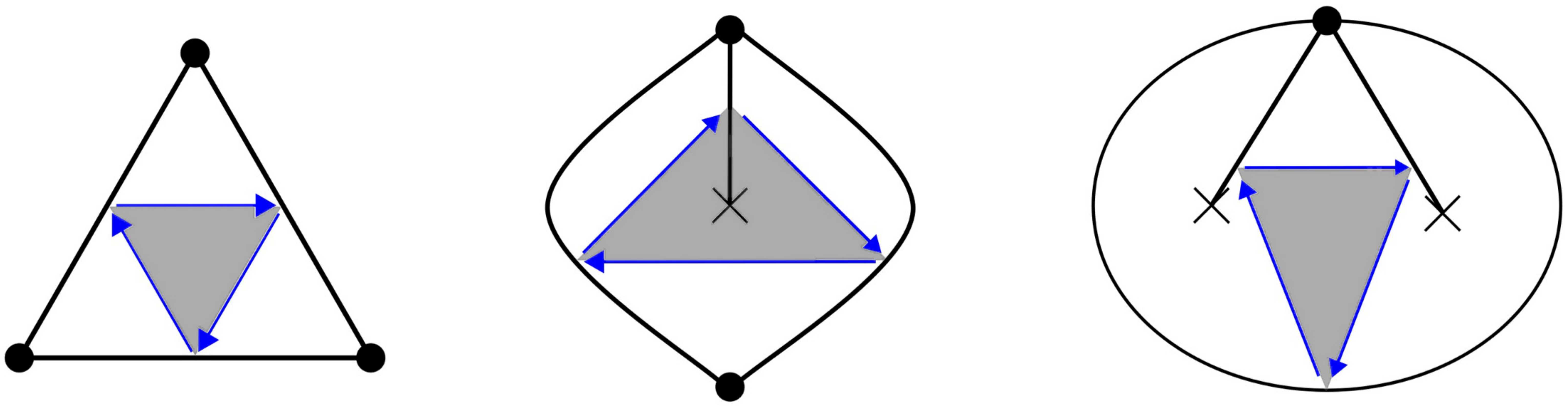}
        \end{figure}

From now on, we shall assume that $\boldsymbol{\Sigma}$ satisfies one of the following conditions:
\begin{equation}\label{DLF_eq:hypotheses-on-surface}
 \partial\Sigma\neq\emptyset \ \text{and} \ \mathbb{M}\subseteq \partial\Sigma \ (\text{so} \ \mathbb{P}=\emptyset), \quad \text{or}\quad
 \partial\Sigma=\emptyset \ \text{and} \ |\mathbb{M}|=1 \ (\text{so} \ \mathbb{P}=\mathbb{M}).
\end{equation}
Then there are three possible types of triangles for any given triangulation $\tau$ of $\boldsymbol{\Sigma}$, namely, those shown in Figure
\ref{DLF_Fig_unpunct_puzzle_pieces}.

\vspace{2mm}

\noindent\textbf{Colored triangulations.} Define a CW-complex $X(\tau)=(X_n(\tau))_{n=0,1,2}$ by:
\begin{align*}
\text{$0$-cells:} &\quad X_0(\tau):=\tau \quad \text{(i.e., each arc in $\tau$ becomes a $0$-cell)}\\
\text{$1$-cells:} &\quad \text{$X_1(\tau):=\{$arrows connecting $0$-cells clockwisely inside triangles$\}$}\\
\text{$2$-cells:} &\quad \text{$X_2(\tau):=\{$$3$-cycles on $X_1(\tau)$ arising from triangles of $\tau$, up to rotation$\}$.}
\end{align*}
Let $C_n(\tau):=\mathbb{F}_2X_n(\tau)$ be the $\mathbb{F}_2$-vector space with basis $X_n(\tau)$, where $\mathbb{F}_2:=\mathbb{Z}/2\mathbb{Z}$, and define a chain complex $C_\bullet(\tau)$ and a cochain complex $C^\bullet(\tau)$ as follows
\begin{align*}
C_\bullet(\tau):\qquad & 0 \longrightarrow C_2(\tau) \overset{\partial_2}{\longrightarrow} C_1(\tau) \overset{\partial_1}{\longrightarrow} C_0(\tau) \longrightarrow 0\\
& \partial_2(\alpha\beta\gamma) :=\alpha+\beta+\gamma, \quad \partial_1(\alpha):=h(\alpha)-t(\alpha)\\
C^\bullet(\tau):=\operatorname{Hom}_{\mathbb{F}_2}(C_\bullet(\tau),\mathbb{F}_2):\qquad & 0 \longleftarrow C^2(\tau) \overset{\partial_2^{\vee}}{\longleftarrow} C^1(\tau) \overset{\partial_1^{\vee}}{\longleftarrow} C^0(\tau) \longleftarrow 0
\end{align*}

Following \cite{DLF_GLF2}, we define a \emph{colored triangulation} of $\boldsymbol{\Sigma}=(\Sigma,\mathbb{M},\mathbb{O})$ to be a pair $(\tau,\xi)$ consisting of a triangulation $\tau$ of $\boldsymbol{\Sigma}$ and a $1$-cocycle $\xi\in \ker(\partial_2^{\vee})\subseteq C^1(\tau)$, i.e. a choice $(\xi_a)_{a\in X_1(\tau)}$ of elements of $\mathbb{F}_2=\{0,1\}$, subject to the condition that for every $2$-cell $\alpha\beta\gamma\in X_2(\tau)$ one must have
$
\xi_\alpha+\xi_\beta+\xi_\gamma=0  \mod 2.
$

\vspace{2mm}

\noindent\textbf{Semilinear clannish algebras.} Let $K$ be any field, 
and suppose we are given:
\begin{itemize}
\item a finite quiver $\widehat{Q}$, not necessarily loop-free;
\item a set $\mathbb{S}\subseteq\widehat{Q}$ of \emph{special loops};
\item a field automorphism $\sigma_a\in\operatorname{Aut}(K)$ for each arrow $a\in\widehat{Q}_1$;
\item a set $Z$ of paths on $\widehat{Q}$ of length at least $2$;
\item a degree-$2$ polynomial $q_s\in K[s;\sigma_s]$ for each $s\in\mathbb{S}$, where $K[s;\sigma_s]$ is the \emph{skew-polynomial ring} in $s$ with coefficients in $K$, skewed by $\sigma_s$.
\end{itemize}
Suppose further that the following conditions are met:
\begin{itemize}
\item for each vertex $k\in \widehat{Q}_0$, at most two arrows of $\widehat{Q}$ end (resp. start) at $k$;
\item for each arrow $a:k\rightarrow j$ not in $\mathbb{S}$, at most one arrow $b$ of $\widehat{Q}$ ends (resp. starts) at $k$ and satisfies $ab\notin Z$ (resp. $ba\notin Z$);
\item no path belonging to $Z$ has a special loop as its first or last arrow.
\end{itemize}
Consider the ring $S:=\times_{k\in \widehat{Q}_0}K$ and the $S$-$S$-bimodule $\bigoplus_{a\in \widehat{Q}_1}K^{\sigma_a}\otimes_KK$, where for any field automophism $\sigma\in\operatorname{Aut}(K)$, we define $K^{\sigma}$ to be the twisted  $K$-$K$-bimodule having $K$ as underlying abelian group, with left $K$-action $z*m:=zm$ and right $K$-action $m*z:=m\sigma(z)$ for $z\in K$ and $m\in K$. Following \cite{DLF_BTCB}, we denote the tensor algebra $K_{\boldsymbol{\sigma}}\widehat{Q}:=T_S(\bigoplus_{a\in \widehat{Q}_1}K^{\sigma_a}\otimes_KK)$,  
and say that the quotient
$$
K_{\boldsymbol{\sigma}}\widehat{Q}/\langle Z\cup \{q_s\mid s\in\mathbb{S}\}\rangle
$$
is a \emph{semilinear clannish algebra}.

As in \cite{DLF_BTLF}, we associate a semilinear clannish algebra to each colored triangulation $(\tau,\xi)$ as follows. Fix one of the following two choices of field $K$:
$$
K :=
\mathbb{C} 
\quad
\text{or}
\quad
K :=
\mathbb{R} 
$$
Set $\widehat{Q}(\tau)$ to be the quiver obtained from $(X_0(\tau),X_1(\tau))$ by adding a loop $s_j$ at each pending arc $j$ of $\tau$. Further, let $\theta:\mathbb{C}\rightarrow\mathbb{C}$ be complex conjugation and define
\begin{align*}
\mathbb{S}(\tau)&:=\{s_j\ | \ j\in\tau \ \text{is pending}\} \qquad \text{i.e., all loops are special};\\
\sigma_a&:=
\begin{cases}
\theta^{\xi_a}|_K & \text{if $a\in X_1(\tau)$}\\
\theta|_K & \text{if $a\in\mathbb{S}(\tau)$}
\end{cases}\qquad \text{for each arrow $a$ of $\widehat{Q}(\tau)$};\\
Z(\tau)&:=\{\alpha\beta\ | \ \exists\gamma\in X_1(\tau) \  \text{such that} \ \alpha\beta\gamma\in X_2(\tau) \ \text{up to rotation of cycles}\},\\
q_{s_j} &:= \begin{cases}
s_j^2-1\in \mathbb{C}[s_j;\theta] & \text{if $K=\mathbb{C}$}\\
s_j^2+1\in \mathbb{R}[s_j;1\hspace{-0.125cm}1_{\mathbb{R}}] & \text{if $K=\mathbb{R}$}
\end{cases}
\qquad \text{for all pending arcs $j$ of $\tau$}. 
\end{align*}
These definitions can be mnemotechnically visualized for the triangles from Figure~\ref{DLF_Fig_unpunct_puzzle_pieces} as follows ($K=\mathbb{C}$ first, $K=\mathbb{R}$ afterwards):
$$\xymatrix{
 & {\mathbb{C}} \ar[dr]|-{\qquad\mathbb{C}^{\theta^{\xi_\beta}}\otimes_\mathbb{C}\mathbb{C}}="b" \\
 {\mathbb{C}} \ar[ur]|-{\mathbb{C}^{\theta^{\xi_\gamma}}\otimes_\mathbb{C}\mathbb{C}\quad}="c" & & {\mathbb{C}} \ar[ll]|-{\mathbb{C}^{\theta^{\xi_\alpha}}\otimes_\mathbb{C}\mathbb{C}}_{ \ }="a"
 \ar @{.}@/^0.65pc/ "c";"b"
 \ar @{.}@/^0.65pc/ "b";"a"
  \ar @{.}@/^0.65pc/ "a";"c"
}\quad
\xymatrix{
 & \mathbb{C} \ar@(ul,ur)^{\mathbb{C}^\theta\otimes_{\mathbb{C}}\mathbb{C}}_{s^2=e \quad s \qquad\quad } \ar[dr]|-{\qquad\mathbb{C}^{\theta^{\xi_\beta}}\otimes_\mathbb{C}\mathbb{C}}="b" \\
 \mathbb{C} \ar[ur]|-{\mathbb{C}^{\theta^{\xi_\gamma}}\otimes_\mathbb{C}\mathbb{C}\quad}="c" & & \mathbb{C} \ar[ll]|-{\mathbb{C}^{\theta^{\xi_\alpha}}\otimes_\mathbb{C}\mathbb{C}}_{ \ }="a"
  \ar @{.}@/^0.65pc/ "c";"b"
 \ar @{.}@/^0.65pc/ "b";"a"
  \ar @{.}@/^0.65pc/ "a";"c"
}
\xymatrix{
\mathbb{C} \ar@(ul,ur)^{\mathbb{C}^\theta\otimes_{\mathbb{C}}\mathbb{C}}_{s_1^2=e_1 \quad s_1 \qquad\quad } \ar[rr]^{\mathbb{C}^{\theta^{\xi_\alpha}}\otimes_\mathbb{C}\mathbb{C}}_{ \ }="d" & & \mathbb{C} \ar@(ul,ur)^{\mathbb{C}^\theta\otimes_{\mathbb{C}}\mathbb{C}}_{\qquad\quad s_2 \quad s_2^2=e_2    } \ar[dl]|-{\quad\mathbb{C}^{\theta^{\xi_\gamma}}\otimes_\mathbb{C}\mathbb{C}}="eta" \\
 & \mathbb{C} \ar[ul]|-{\mathbb{C}^{\theta^{\xi_{\beta}}}\otimes_\mathbb{C}\mathbb{C}\qquad}="eps" 
 \ar @{.}@/^0.65pc/ "eps";"d"
 \ar @{.}@/^0.65pc/ "d";"eta"
  \ar @{.}@/^0.65pc/ "eta";"eps"
}
$$
$$\xymatrix{
 & {\mathbb{R}} \ar[dr]|-{\qquad\mathbb{R}^{{}}\otimes_\mathbb{R}\mathbb{R}}="b" \\
 {\mathbb{R}} \ar[ur]|-{\mathbb{R}^{{}}\otimes_\mathbb{R}\mathbb{R}\qquad}="c" & & {\mathbb{R}} \ar[ll]|-{\mathbb{R}^{{}}\otimes_\mathbb{R}\mathbb{R}}_{ \ }="a"
 \ar @{.}@/^0.65pc/ "c";"b"
 \ar @{.}@/^0.65pc/ "b";"a"
  \ar @{.}@/^0.65pc/ "a";"c"
}\quad
\xymatrix{
 & \mathbb{R} \ar@(ul,ur)^{\mathbb{R}\otimes_{\mathbb{R}}\mathbb{R}}_{s^2=-e \quad s \qquad\quad \ } \ar[dr]|-{\qquad\mathbb{R}\otimes_\mathbb{R}\mathbb{R}}="b" \\
 \mathbb{R} \ar[ur]|-{\mathbb{R}\otimes_\mathbb{R}\mathbb{R}\qquad}="c" & & \mathbb{R} \ar[ll]|-{\mathbb{R}\otimes_\mathbb{R}\mathbb{R}}_{ \ }="a"
  \ar @{.}@/^0.65pc/ "c";"b"
 \ar @{.}@/^0.65pc/ "b";"a"
  \ar @{.}@/^0.65pc/ "a";"c"
}
\xymatrix{
\mathbb{R} \ar@(ul,ur)^{\mathbb{R}\otimes_{\mathbb{R}}\mathbb{R}}_{s_1^2=-e_1 \quad s_1 \qquad\qquad } \ar[rr]^{\mathbb{R}\otimes_\mathbb{R}\mathbb{R}}_{ \ }="d" & & \mathbb{R} \ar@(ul,ur)^{\mathbb{R}\otimes_{\mathbb{R}}\mathbb{R}}_{\qquad\qquad s_2 \quad s_2^2=-e_2    } \ar[dl]|-{\quad\mathbb{R}\otimes_\mathbb{R}\mathbb{R}}="eta" \\
 & \mathbb{R} \ar[ul]|-{\mathbb{R}\otimes_\mathbb{R}\mathbb{R}\quad}="eps" 
 \ar @{.}@/^0.65pc/ "eps";"d"
 \ar @{.}@/^0.65pc/ "d";"eta"
  \ar @{.}@/^0.65pc/ "eta";"eps"
}
$$

We show in \cite{DLF_BTLF} that $K_{\boldsymbol{\sigma}}\widehat{Q}(\tau)/I(\tau,\xi)$, where $I(\tau,\xi):=\langle Z\cup \{q_{s_j}\ | \ s_j\in\mathbb{S}(\tau)\}\rangle$, is a semilinear clannish algebra.

\vspace{2mm}

\noindent \textbf{Jacobian algebras.} Let  $(\tau,\xi)$ be a colored triangulation of $\boldsymbol{\Sigma}$. 
Pick one of the following two assignments $\mathbf{F}=(F_k)_{k\in \tau}$ of fields $F_k$ for $k\in\tau$:
\begin{equation*}
F_k :=\begin{cases}
\mathbb{R} & \text{for all $k$ pending};\\
\mathbb{C} & \text{for all $k$ non-pending};
\end{cases}
\
\text{or}
\
F_k :=\begin{cases}
\mathbb{C} & \text{for all $k$ pending};\\
\mathbb{R} & \text{for all $k$ non-pending}.
\end{cases}
\end{equation*}
We shall say that the assignment on the left \emph{is $B$-like}, and to that the one on the right \emph{is $C$-like}.
For each arrow $a\in X_1(\tau)$, $a:k\rightarrow j$, set $g(\tau,\xi)_a:=\theta^{\xi_a}|_{F_j\cap F_k}$ and
\begin{equation*}
A(\tau,\xi)_a :=\begin{cases}
(F_j\otimes_{\mathbb{R}}F_k)^2=(\mathbb{R}\otimes_{\mathbb{R}}\mathbb{R})^2 & \text{if $k,j$ are pending and $\mathbf{F}$ is $B$-like};\\
F_j\otimes_{\mathbb{R}}F_k=\mathbb{C}\otimes_{\mathbb{R}}\mathbb{C} & \text{if $k,j$ are pending and $\mathbf{F}$ is $C$-like};\\
F_j^{g(\tau,\xi)_a}\otimes_{F_j\cap F_k}F_k & \text{if  at least one of $k,j$ is non-pending}.
\end{cases}
\end{equation*}
These assignments of fields and bimodules can be mnemotechnically visualized for the triangles from Figure \ref{DLF_Fig_unpunct_puzzle_pieces} as follows ($B$-like assignment first, $C$-like afterwards)
$$\xymatrix{
 & {\mathbb{C}} \ar[dr]|-{\qquad\mathbb{C}^{\theta^{\xi_\beta}}\otimes_\mathbb{C}\mathbb{C}}="b" \\
 {\mathbb{C}} \ar[ur]|-{\mathbb{C}^{\theta^{\xi_\gamma}}\otimes_\mathbb{C}\mathbb{C}\quad}="c" & & {\mathbb{C}} \ar[ll]|-{\mathbb{C}^{\theta^{\xi_\alpha}}\otimes_\mathbb{C}\mathbb{C}}_{ \ }="a"
}\quad
\xymatrix{
 & \mathbb{R}  \ar[dr]|-{\quad\mathbb{C}\otimes_\mathbb{R}\mathbb{R}}="b" \\
 \mathbb{C} \ar[ur]|-{\mathbb{R}\otimes_\mathbb{R}\mathbb{C}\ }="c" & & \mathbb{C} \ar[ll]|-{\mathbb{C}^{\theta^{\xi_\alpha}}\otimes_\mathbb{C}\mathbb{C}}_{ \ }="a"
}
\xymatrix{
\mathbb{R}  \ar[rr]|-{(\mathbb{R}\otimes_\mathbb{R}\mathbb{R})^2}_{ \ }="d" & & \mathbb{R}  \ar[dl]|-{\quad\mathbb{C}\otimes_\mathbb{R}\mathbb{R}}="eta" \\
 & \mathbb{C} \ar[ul]|-{\mathbb{R}\otimes_\mathbb{R}\mathbb{C}\quad}="eps" 
}
$$
$$\xymatrix{
 & {\mathbb{R}} \ar[dr]|-{\qquad\mathbb{R}^{{}}\otimes_\mathbb{R}\mathbb{R}}="b" \\
 {\mathbb{R}} \ar[ur]|-{\mathbb{R}^{{}}\otimes_\mathbb{R}\mathbb{R}\qquad}="c" & & {\mathbb{R}} \ar[ll]|-{\mathbb{R}^{{}}\otimes_\mathbb{R}\mathbb{R}}_{ \ }="a"
}\quad
\xymatrix{
 & \mathbb{C} \ar[dr]|-{\qquad\mathbb{R}\otimes_\mathbb{C}\mathbb{C}}="b" \\
 \mathbb{R} \ar[ur]|-{\mathbb{C}\otimes_\mathbb{R}\mathbb{R}\qquad}="c" & & \mathbb{R} \ar[ll]|-{\mathbb{R}\otimes_\mathbb{R}\mathbb{R}}_{ \ }="a"
}
\xymatrix{
\mathbb{C}  \ar[rr]|-{\mathbb{C}\otimes_\mathbb{R}\mathbb{C}}_{ \ }="d" & & \mathbb{C}  \ar[dl]|-{\quad\mathbb{R}\otimes_\mathbb{R}\mathbb{C}}="eta" \\
 & \mathbb{R} \ar[ul]|-{\mathbb{C}\otimes_\mathbb{R}\mathbb{R}\quad}="eps" 
}
$$

Let $R:=\times_{k\in\tau}F_k$, $A(\tau,\xi):=\bigoplus_{a\in X_1(\tau)}A(\tau,\xi)_a$. Thus, $A(\tau,\xi)$ is an $R$-$R$-bimodule, so one can form the (complete) tensor ring of $A(\tau,\xi)$ over $R$. Following \cite{DLF_GLF1,DLF_GLF2}, we can define an ``obvious'' (super-)potential $W(\tau,\xi)\in T_R(A(\tau,\xi))$ as the sum of ``obvious'' degree-$3$ cycles in $T_R(A(\tau,\xi))$, and take the \emph{cyclic derivatives} of $W(\tau,\xi)$ to define the \emph{Jacobian algebra} $\mathcal{P}(A(\tau,\xi),W(\tau,\xi))$. 
Roughly, for each arrow $a$ and cycle $c$, the cyclic derivative $\partial_a(c)$ is defined to be the $g(\tau,\xi)_a^{-1}$-linear part of the usual sum of paths obtained by deleting each occurrence of $a$ in $c$ (with the reordering $yx$ if $c=xay$),  see \cite{DLF_GLF1}. Since the $F_k$ are not necessarily all the same field, the notion of path has to be enhanced, so both $a$ and $ia$ (resp. $ai$) may be paths whenever $a$ is an arrow with $F_{h(a)}=\mathbb{C}$ (resp. $F_{t(a)}=\mathbb{C}$). This way, e.g., $ab$ and $aib$ are distinct paths if $t(a)=h(b)$, $F_{h(a)}=\mathbb{R}$, $F_{t(a)}=\mathbb{C}$, $F_{t(b)}=\mathbb{R}$. On the other hand, if $t(a)=h(b)$, $F_{h(a)}=\mathbb{C}$, $F_{t(a)}=\mathbb{R}$, $F_{t(b)}=\mathbb{C}$, then
\begin{align*}
ab&=\frac{1}{2}(ab-iabi)+\frac{1}{2}(ab+iabi), & \\
 (ab-iabi)i &= i(ab-iabi), & (ab+iabi)i=-i(ab+iabi),
\end{align*}
i.e., $\frac{1}{2}(ab-iabi)$ and $\frac{1}{2}(ab+iabi)$ are the $1\hspace{-0.125cm}1$-linear part and the $\theta$-linear part of $ab$.

\vspace{2mm}

\noindent\textbf{Main result.} Essential to the proof of our main result are the simple observations that
$
\mathbb{C}  \cong \mathbb{R}[s]/\langle s^2+1\rangle=\mathbb{R}[s;1\hspace{-0.125cm}1_{\mathbb{R}}]/\langle s^2+1\rangle$ and $
\mathbb{R}  \simeq_{\operatorname{Morita}} \mathbb{C}[s;\mathbb{\theta}]/\langle s^2-1\rangle.
$

\begin{theorem}\cite{DLF_BTLF}\label{DLF:main-result} If $\boldsymbol{\Sigma}$ is a surface with marked points and orbifold points satisfying  \eqref{DLF_eq:hypotheses-on-surface}, then for any colored triangulation $(\tau,\xi)$ of $\boldsymbol{\Sigma}$, the Jacobian algebra $\mathcal{P}(A(\tau,\xi),W(\tau,\xi))$ and the semilinear clannish algebra $K_{\boldsymbol{\sigma}}\widehat{Q}(\tau,\xi)/I(\tau,\xi)$ are Morita-equivalent ($K=\mathbb{C}$ in the $B$-like situation, $K=\mathbb{R}$ in the $C$-like situation).
\end{theorem}

\end{document}